\documentclass[11pt,a4paper]{amsart}
\usepackage[margin=30mm]{geometry}
\usepackage[T1]{fontenc}
\usepackage{amsmath,amssymb,amsthm,mathtools}
\usepackage{enumitem,booktabs,array,tabularx,aliascnt}
\usepackage{xcolor}
\usepackage{hyperref}
\usepackage[nameinlink,capitalize,noabbrev]{cleveref}
\hypersetup{colorlinks=true,linkcolor=blue,citecolor=blue,urlcolor=blue,pdftitle={On the h-polynomials of cyclotomic standard graded Cohen-Macaulay rings},pdfauthor={Akihiro Higashitani and Takumi Yanagida}}
\setlist[enumerate]{leftmargin=2.2em}
\setlist[itemize]{leftmargin=2em}
\allowdisplaybreaks
\newtheorem{theorem}{Theorem}[section]
\newaliascnt{proposition}{theorem}
\newtheorem{proposition}[proposition]{Proposition}
\aliascntresetthe{proposition}
\newaliascnt{lemma}{theorem}
\newtheorem{lemma}[lemma]{Lemma}
\aliascntresetthe{lemma}
\newaliascnt{corollary}{theorem}

\aliascntresetthe{corollary}
\theoremstyle{definition}
\newaliascnt{definition}{theorem}
\newtheorem{definition}[definition]{Definition}
\aliascntresetthe{definition}
\newaliascnt{question}{theorem}
\newtheorem{question}[question]{Question}
\aliascntresetthe{question}
\theoremstyle{remark}
\newaliascnt{remark}{theorem}
\newtheorem{remark}[remark]{Remark}
\aliascntresetthe{remark}
\crefname{theorem}{Theorem}{Theorems}
\Crefname{theorem}{Theorem}{Theorems}
\crefname{proposition}{Proposition}{Propositions}
\Crefname{proposition}{Proposition}{Propositions}
\crefname{lemma}{Lemma}{Lemmas}
\Crefname{lemma}{Lemma}{Lemmas}
\crefname{corollary}{Corollary}{Corollaries}
\Crefname{corollary}{Corollary}{Corollaries}
\crefname{definition}{Definition}{Definitions}
\Crefname{definition}{Definition}{Definitions}
\crefname{question}{Question}{Questions}
\Crefname{question}{Question}{Questions}
\crefname{remark}{Remark}{Remarks}
\Crefname{remark}{Remark}{Remarks}

\newcommand{\kk}{\Bbbk}
\newcommand{\ZZ}{\mathbb Z}
\newcommand{\QQ}{\mathbb Q}
\DeclareMathAlphabet{\cal}{OMS}{cmsy}{m}{n}
\newcommand{\Nc}{{\cal N}}
\DeclareMathOperator{\Hilb}{Hilb}
\newcommand{\ph}{\varphi}

\title[Cyclotomic $h$-polynomials of Cohen--Macaulay rings]{On the $h$-polynomials of cyclotomic standard graded Cohen--Macaulay rings}
\author{Akihiro Higashitani}
\address{Department of Pure and Applied Mathematics, Graduate School of Information Science and Technology, Osaka University, Suita, Osaka 565-0871, Japan}
\email{higashitani@ist.osaka-u.ac.jp}
\author{Takumi Yanagida}
\address{Department of Pure and Applied Mathematics, Graduate School of Information Science and Technology, Osaka University, Suita, Osaka 565-0871, Japan}
\email{u465270i@ecs.osaka-u.ac.jp}

\begin{document}
\begin{abstract}
We say that a standard graded algebra is cyclotomic if all roots of its $h$-polynomial lie on the complex unit circle. We study which products of cyclotomic polynomials can occur as $h$-polynomials of standard graded Cohen--Macaulay algebras without being of complete intersection type. For each fixed positive integer $a$, we reduce the nonexistence problem for multiplicities $ap$ over all sufficiently large primes $p$ to finitely many cases. Combining this reduction with explicit constructions, we completely characterize the positive integers $k$ for which such a polynomial of multiplicity $8k$ exists: precisely those satisfying
\[k\notin\{1,2,3,6\}\cup\{p:p\ge7\text{ is prime}\}\cup\{3p:p\ge17\text{ is prime}\}.\]
\end{abstract}
\maketitle

\section{Introduction}
Throughout this paper, let $\kk$ be a field and let all rings be commutative. A \emph{standard graded $\kk$-algebra} is a graded algebra $R=\bigoplus_{n\ge0}R_n$ with $R_0=\kk$ that is finitely generated by $R_1$. Its Hilbert function is $H_R(n)=\dim_{\kk}R_n$, and its Hilbert series has the form
\[\Hilb(R,t)=\sum_{n\ge0}H_R(n)t^n=\frac{h_R(t)}{(1-t)^d},\]
where $d=\dim R$ and $h_R(t)\in\ZZ[t]$ satisfies $h_R(0)=1$ and $h_R(1)>0$. The numerator $h_R(t)$ is called the \emph{$h$-polynomial} of $R$, and its coefficient sequence is the \emph{$h$-vector}. The value $h_R(1)$ is the multiplicity of $R$; see \cite[Section~4.1]{BrunsHerzog}.

A standard graded algebra is said to be \emph{cyclotomic} if all roots of its $h$-polynomial lie on the complex unit circle. By Kronecker's theorem, these $h$-polynomials are products of cyclotomic polynomials.
For $m\ge1$, set \[\Psi_m(t)=1+t+\cdots+t^{m-1}.\]
We say that a polynomial is \emph{of type CI} if it is a finite product of polynomials $\Psi_m(t)$. This terminology reflects the fact that the $h$-polynomial of a standard graded complete intersection has this form.

The study of cyclotomic Hilbert series connects the factorization of $h$-polynomials with structural properties of graded algebras. Borz\`i and D'Al\`i proved that a standard graded Koszul algebra is cyclotomic if and only if it is a complete intersection \cite[Theorem~3.9]{BorziDAli}. They also characterized the case of an irreducible cyclotomic $h$-polynomial, showing that the algebra is a hypersurface \cite[Theorem~4.3]{BorziDAli}, and observed that a cyclotomic Cohen--Macaulay graded domain is Gorenstein \cite[Proposition~2.7]{BorziDAli}. These results indicate that cyclotomicity can impose strong algebraic restrictions, even though the condition alone does not characterize complete intersections.

Related questions appear in other settings. Ciolan, Garc\'ia-S\'anchez, and Moree \cite{CiolanGarciaMoree} conjectured that a numerical semigroup is cyclotomic precisely when it is a complete intersection. In Ehrhart theory, Braun and Liu \cite{BraunLiu} investigated lattice simplices with $h^*$-polynomials factoring into products of geometric series. These connections motivate the study of which products of cyclotomic polynomials occur as $h$-polynomials of standard graded algebras. In contrast to these results for domains and Koszul algebras, we impose neither assumption and study the type-CI property at the level of $h$-polynomials.

Higashitani and Ueyama \cite[Theorem~1.3]{HigashitaniUeyama} proved that the $h$-polynomial of a cyclotomic standard graded algebra is of type CI whenever its multiplicity is $1$, $4$, $6$, or a prime. Conversely, for every composite integer $n\ge8$, they constructed a cyclotomic standard graded algebra whose $h$-polynomial is non-type-CI and has multiplicity $n$. This leads to the following question.
\begin{question}
Which multiplicities occur for non-type-CI $h$-polynomials of cyclotomic standard graded Cohen--Macaulay rings?
\end{question}

The main purpose of this paper is to establish a Cohen--Macaulay analogue of \cite[Theorem~1.3]{HigashitaniUeyama}. The additional Cohen--Macaulay assumption changes the answer substantially and makes the problem more delicate. For example, no product of cyclotomic polynomials of multiplicity $8,12,16,20,24$, or $28$ which is non-type-CI admits a Cohen--Macaulay realization (\cref{lem:finite-checks}\textup{(i)}), while all these multiplicities occur without the Cohen--Macaulay assumption by \cite[Theorem~1.3(2)]{HigashitaniUeyama}.

We call a polynomial $f(t)$ \emph{CM-realizable} if it is the $h$-polynomial of a standard graded Cohen--Macaulay algebra. Equivalently, its coefficient sequence is an $M$-sequence, or it is the $h$-polynomial of a Cohen--Macaulay simplicial complex; see \cref{thm:macaulay}.

Our nonexistence results rely on reducing the problem for infinitely many prime parameters to a finite set of polynomials.
For a positive integer $a$, let \[p(a)=\min\{p:p\text{ is prime and }p>a\}.\]
\begin{theorem}\label{thm:prime-reduction}
Let $a$ be a positive integer and let $p\ge p(a)$ be prime. Suppose that there exists a product $f(t)$ of cyclotomic polynomials which is CM-realizable and non-type-CI and satisfies $f(1)=ap$. Then there exist a positive divisor $\ell$ of $a$ and a product $C(t)$ of cyclotomic polynomials satisfying
\[C(1)=\frac{a}{\ell},\qquad \deg C(t)\le\frac{a}{\ell}-1\]
such that $\widetilde f(t)=C(t)\Psi_{\ell p(a)}(t)$ is CM-realizable and non-type-CI.
\end{theorem}
For each fixed $a$, only finitely many pairs $(\ell,C(t))$ satisfy the displayed conditions. Thus, if none of these pairs produces a polynomial $\widetilde f(t)$ which is CM-realizable and non-type-CI, then no product of cyclotomic polynomials with both properties has multiplicity $ap$ for any prime $p\ge p(a)$. The proof of \cref{thm:prime-reduction} is given in \cref{sec:prime-proof}.

Combining this reduction with explicit constructions, we obtain a complete classification at multiplicities $8k$. Set
\begin{equation}\label{eq:exceptional-set}
\Nc=\{1,2,3,6\}\cup\{p:p\ge7\text{ is prime}\}\cup\{3p:p\ge17\text{ is prime}\}.
\end{equation}
\begin{theorem}\label{thm:8k-classification}
Let $k$ be a positive integer. Then there exists a product $f(t)$ of cyclotomic polynomials which is CM-realizable and non-type-CI and satisfies $f(1)=8k$ if and only if $k\notin\Nc$.
\end{theorem}
The existence proof uses constructions for every multiplicity $32k$, for $16k$ with $k$ odd and $k\ne1,3$, and for $8pq$ with $p,q\ge5$ odd integers, together with the three explicit polynomials corresponding to $k=21,33,39$. We apply \cref{thm:prime-reduction} with $a=8,24$ for the nonexistence proof and treat $24p$ for $p=17,19,23$ separately.

The paper is organized as follows. In \cref{sec:preliminaries}, we recall cyclotomic polynomials, Macaulay's theorem, and the operations on simplicial complexes used in the constructions. In Section~\ref{sec:prime-proof}, we prove Theorem~\ref{thm:prime-reduction}, and in \cref{sec:8k}, we establish the classification at multiplicities $8k$. Section~\ref{sec:further} contains partial results for $2k$ and $4k$, the smallest odd multiplicity, and a comparison with realizations by Gorenstein Stanley--Reisner rings. The finite computations used in these results are collected in Appendix~\ref{sec:computation}.

\medskip
\noindent\textbf{Acknowledgments.} The first author was partially supported by JSPS KAKENHI Grant Numbers JP24K00521 and JP24K00534.

\bigskip

\section{Preliminaries}\label{sec:preliminaries}
In this section, we recall the basic facts about cyclotomic polynomials, Macaulay's theorem, and simplicial complexes that will be used throughout the paper.

\subsection{Cyclotomic polynomials and polynomials of type CI}\label{sec:cyclotomic}
For $m\ge2$, the $m$-th cyclotomic polynomial is
\[\Phi_m(t)=\prod_{\substack{1\le j<m\\\gcd(j,m)=1}}\bigl(t-\exp(2\pi j\sqrt{-1}/m)\bigr).\]
It is monic and irreducible in $\QQ[t]$, has integer coefficients and constant term $1$, and is palindromic of degree $\ph(m)$, where $\ph$ is Euler's totient function. Here, a polynomial $f(t)=\sum_{i=0}^d c_it^i$ is called \emph{palindromic} if $c_i=c_{d-i}$ for all $i$. We shall use
\[\Psi_m(t)=\prod_{\substack{r\mid m\\r>1}}\Phi_r(t),\qquad
\Phi_m(1)=\begin{cases}p,&m=p^e\text{ for a prime }p\text{ and }e\ge1,\\1,&m\text{ is not a prime power}.\end{cases}\]
In particular, $\Psi_p(t)=\Phi_p(t)$ holds for every prime $p$. For more details, see, e.g., \cite{Sanna}.

A \emph{product of cyclotomic polynomials} will always mean $f(t)=\prod_{m>1}\Phi_m(t)^{e_m}$, where $e_m\in\ZZ_{\ge0}$ and all but finitely many $e_m$ vanish. Such a product is monic, has constant term $1$, and is palindromic.

\begin{definition}\label{def:type-ci}
A product $f(t)$ of cyclotomic polynomials is called \emph{of type CI} if it admits a factorization $f(t)=\prod_{m>1}\Psi_m(t)^{c_m}$ with nonnegative integers $c_m$, all but finitely many of which are zero. Otherwise, we call it \emph{non-type-CI}.
\end{definition}

\begin{lemma}\label{lem:mobius-ci}
Let $f(t)=\prod_{m>1}\Phi_m(t)^{e_m}$ and set \[b_n(f(t))=\sum_{r\ge1}\mu(r)e_{nr} \;\text{ for }\; n>1,\] where $\mu$ is the M\"obius function.
Then $f(t)$ is of type CI if and only if $b_n(f(t))\ge0$ for all $n>1$. In that case, $f(t)=\prod_{n>1}\Psi_n(t)^{b_n(f(t))}$.
\end{lemma}
\begin{proof}
The identity $f(t)=\prod_{n>1}\Psi_n(t)^{c_n}$ holds if and only if $e_m=\sum_{r\ge1}c_{mr}$ for all $m>1$, by the cyclotomic factorization of $\Psi_n(t)$. M\"obius inversion transforms these equalities into $c_n=\sum_{r\ge1}\mu(r)e_{nr}=b_n(f(t))$. Hence, a type-CI factorization exists precisely when all these integers are nonnegative.
\end{proof}
In particular, multiplication by $\Psi_m(t)$ increases $b_m(f(t))$ by $1$ and leaves every other $b_n(f(t))$ unchanged.

\subsection{Macaulay's inequalities and CM realizability}\label{sec:M}
For integers $c>0$ and $i\ge1$, the $i$-th \emph{Macaulay representation} of $c$ is the unique expression
\[c=\binom{c_i}{i}+\binom{c_{i-1}}{i-1}+\cdots+\binom{c_j}{j},\qquad c_i>c_{i-1}>\cdots>c_j\ge j\ge1.\]
Set \[c^{\langle i\rangle}=\binom{c_i+1}{i+1}+\binom{c_{i-1}+1}{i}+\cdots+\binom{c_j+1}{j+1}\;\;\text{ and }\;\; 0^{\langle i\rangle}=0.\]

\begin{definition}
A sequence $H=(H(0),H(1),\ldots)$ of nonnegative integers is an \emph{$M$-sequence} if $H(0)=1$ and $H(i+1)\le H(i)^{\langle i\rangle}$ for every $i\ge1$. A finite sequence is called an $M$-sequence if its extension by zeros satisfies this condition.
\end{definition}

\begin{theorem}[{cf. \cite[Theorems~4.2.10 and~5.1.15]{BrunsHerzog}}]\label{thm:macaulay}
A sequence $H=(H(0),H(1),\ldots)$ of nonnegative integers is the Hilbert function of a standard graded $\kk$-algebra if and only if it is an $M$-sequence. Moreover, a finite sequence of nonnegative integers with first entry $1$ is the $h$-vector of a Cohen--Macaulay simplicial complex over $\kk$ if and only if it is an $M$-sequence.
\end{theorem}
Consequently, a polynomial is CM-realizable if and only if its coefficient sequence is a finite $M$-sequence. In fact, the $h$-vector of a standard graded Cohen--Macaulay ring is the Hilbert function of an Artinian reduction (after extending the field if necessary); conversely, a finite $M$-sequence is itself the Hilbert function, and hence the $h$-vector, of an Artinian standard graded ring. The characterization in terms of simplicial complexes provides a realization over any prescribed field.

We shall repeatedly use the identity
\begin{equation}\label{eq:small-macaulay}
c^{\langle i\rangle}=c\;\;\text{ for }\;\; 0\le c\le i.
\end{equation}
For $c>0$, this follows from $c=\binom{i}{i}+\cdots+\binom{i-c+1}{i-c+1}$; the case $c=0$ follows from the definition. Thus, an $M$-sequence cannot increase after a position at which its value does not exceed its index. Once an entry is zero, every subsequent entry is zero. We also use $c^{\langle i\rangle}\ge c$. In particular, a nonincreasing step always satisfies Macaulay's inequality.

\begin{lemma}[{\cite[Proposition~3.2]{HigashitaniUeyama}}]\label{lem:degree-bound}
If a product $f(t)$ of cyclotomic polynomials is the $h$-polynomial of a standard graded algebra, then $\deg f(t)\le f(1)-1$. In particular, this is true whenever $f(t)$ is CM-realizable.
\end{lemma}
For a polynomial which is CM-realizable, this bound also follows immediately from the absence of internal zeros in its coefficient sequence.

The next estimate will be used to verify an infinite family of coefficient sequences in Section~\ref{sec:8k}.
\begin{lemma}\label{lem:macaulay-increment}
Let $i\ge1$, let $h_i,h_{i+1}$ be nonnegative integers, and set $d=h_{i+1}-h_i$. If either $d\le0$, or $1\le d\le i$ and $h_i\ge di-\binom{d-1}{2}+1$, then $h_{i+1}\le h_i^{\langle i\rangle}$.
\end{lemma}
\begin{proof}
The case $d\le0$ follows from $h_i^{\langle i\rangle}\ge h_i$. Suppose $1\le d\le i$ and set $T_{i,d}=di-\binom{d-1}{2}+1$. Its Macaulay representation is
\[T_{i,d}=\binom{i+1}{i}+\binom{i}{i-1}+\cdots+\binom{i-d+2}{i-d+1},\]
so $T_{i,d}^{\langle i\rangle}=T_{i,d}+d$. By \cite[Lemma~4.2.13(b)]{BrunsHerzog}, $(a+1)^{\langle i\rangle}\ge a^{\langle i\rangle}+1$ holds for every positive integer $a$. Hence, $a\mapsto a^{\langle i\rangle}-a$ is nondecreasing.
Therefore, \[h_i^{\langle i\rangle}-h_i\ge T_{i,d}^{\langle i\rangle}-T_{i,d}=d,\]
which is the required inequality.
\end{proof}

\subsection{Simplicial complexes and Stanley--Reisner rings}
We recall the notation needed for the constructions; see, e.g., \cite[Section~5.1]{BrunsHerzog}. A \emph{simplicial complex} $\Delta$ on a finite vertex set $V$ is a collection of subsets of $V$ containing the empty set and every singleton, and closed under taking subsets. Its dimension is defined by $\max\{|F|-1:F\in\Delta\}$. For $V=[n]=\{1,\ldots,n\}$, its Stanley--Reisner ring is
\[\kk[\Delta]=\kk[x_1,\ldots,x_n]/I_\Delta,\;\;\text{ where }\;\; I_\Delta=(\,\textstyle\prod_{i\in F}x_i:F\subseteq[n],\ F\notin\Delta\,).\]
The ideal $I_\Delta$ is called the \emph{Stanley--Reisner ideal}. We say that $\Delta$ is Cohen--Macaulay over $\kk$ when $\kk[\Delta]$ is Cohen--Macaulay.
If $\dim\Delta=d-1$, then its $h$-polynomial is \[h_\Delta(t)=h_{\kk[\Delta]}(t)=\sum_{i=0}^d f_{i-1}(\Delta)t^i(1-t)^{d-i},\]
where $f_{i-1}(\Delta)$ is the number of faces of cardinality $i$ and $f_{-1}(\Delta)=1$.

For simplicial complexes $\Delta$ and $\Gamma$ on disjoint vertex sets, their \emph{join} is $\Delta*\Gamma=\{F\cup G:F\in\Delta,\ G\in\Gamma\}$. It satisfies $h_{\Delta*\Gamma}(t)=h_\Delta(t)h_\Gamma(t)$ and is Cohen--Macaulay if both factors are. The boundary complex $\partial\Delta_d$ of a $d$-simplex has $h$-polynomial $\Psi_{d+1}(t)$. It follows that multiplication by any $\Psi_m(t)$ preserves CM realizability.

\begin{lemma}\label{lem:missing-divisor}
Let $f(t)$ be a product of cyclotomic polynomials. If there exist integers $r,m\ge2$ such that $r\mid m$, $\Phi_m(t)\mid f(t)$, and $\Phi_r(t)\nmid f(t)$, then $f(t)$ is non-type-CI.
\end{lemma}
\begin{proof}
If $f(t)$ admits a type-CI factorization, then a factor $\Psi_n(t)$ divisible by $\Phi_m(t)$ must satisfy $m\mid n$. Since $r \mid m$, the same factor contains $\Phi_r(t)$, a contradiction.
\end{proof}

\bigskip

\section{Proof of the main theorem}\label{sec:prime-proof}
The goal of this section is to prove Theorem~\ref{thm:prime-reduction}. We begin with two elementary reductions.
\begin{lemma}\label{lem:p-primary-factor}
Let $a$ be a positive integer and let $p>a$ be prime. If $f(t)$ is a product of cyclotomic polynomials which is CM-realizable and satisfies $f(1)=ap$, then $f(t)=A(t)\Phi_p(t)$ for a product $A(t)$ of cyclotomic polynomials with $A(1)=a$.
\end{lemma}
\begin{proof}
Since $p\nmid a$, exactly one cyclotomic factor of $f(t)$ has index a power of $p$, counted with multiplicity. Write it as $\Phi_{p^e}(t)$. If $e\ge2$, then $\deg f(t)\ge p(p-1)\ge ap$, a contradiction to \cref{lem:degree-bound}. Thus, $e=1$, and we obtain $A(t)$ by removing this factor.
\end{proof}

\begin{lemma}\label{lem:small-degree-A}
Let $m\ge2$ be an integer and let $A(t)$ be a product of cyclotomic polynomials with $A(1)=a$ and $d=\deg A(t)<m$. If $A(t)\Psi_m(t)$ is CM-realizable, then $d\le a-1$.
\end{lemma}
\begin{proof}
If $d=0$, then $A(t)=1$; if $d=1$, then $A(t)=1+t$. The assertion follows in both cases. Suppose $d\ge2$. Since $d<m$, the coefficients of $t^{d-1}$ and $t^d$ in $A(t)\Psi_m(t)$ are $a-1$ and $a$, respectively. If $d\ge a$, then Macaulay's inequality at degree $d-1$ would imply $a\le(a-1)^{\langle d-1\rangle}=a-1$, a contradiction.
\end{proof}

When we consider $f(t)=A(t)\Phi_p(t)$, we write
\begin{equation}\label{eq:A-f-notation}
A(t)=\sum_{j=0}^d a_jt^j,\qquad f(t)=\sum_{i=0}^R h_it^i,\qquad R=d+p-1,
\end{equation}
and set coefficients outside the indicated ranges equal to zero. In particular, $h_i-h_{i-1}=a_i-a_{i-p}$ for $i\ge1$.

\begin{lemma}\label{lem:constant-interval}
Let $a$ be a positive integer, let $p>a$ be prime, and let $A(t)$ be a product of cyclotomic polynomials with $A(1)=a$ and $d=\deg A(t)\ge p$. If $f(t)=A(t)\Phi_p(t)$ is CM-realizable, then the least index $s$ with $h_s\le s$ satisfies
\[1\le s<p,\qquad h_i>i\ (0\le i<s),\qquad h_{s-1}=h_s=\cdots=h_{R-s+1}=s.\]
\end{lemma}
\begin{proof}
Since $R=d+p-1 \ge2p-1$, the index sets $\{0,\ldots,p-1\}$ and $\{R-p+1,\ldots,R\}$ are disjoint. If $h_i>i$ for every $0\le i<p$, then palindromicity implies
\[ap=\sum_{i=0}^R h_i\ge2\sum_{i=0}^{p-1}(i+1)=p(p+1)>ap,\]
a contradiction. Thus, $s<p$, and $h_0=1$ ensures $s\ge1$.

Starting from $h_s\le s$, Macaulay's inequalities force all subsequent coefficients to be nonincreasing. In particular, $h_s\ge h_{s+1}\ge\cdots\ge h_{R-s}=h_s$, where the last equality follows from palindromicity. Moreover, minimality of $s$ implies $h_{s-1}\ge s$, while palindromicity implies $h_{R-s+1}=h_{s-1}$.
Therefore, \[h_{s-1} \geq s \geq h_s \geq \cdots \geq h_{R-s} = h_s \geq h_{R-s+1}=h_{s-1},\]
so every inequality is an equality, as required.
\end{proof}
\begin{lemma}\label{lem:residue-decomposition}
Let $p$ be prime and let $A(t)=\sum_{j=0}^d a_jt^j$ be a product of cyclotomic polynomials with $d\ge p$ and $A(1)=a$. Write $A(t)\Phi_p(t)=\sum_{i=0}^R h_it^i$, where $R=d+p-1$, and suppose that $h_s=\cdots=h_{R-s}$ for an integer $1\le s<p$. Write $d=qp+r$ with $q\ge1$ and $0\le r<p$, and set
\[\alpha=p+r-s,\qquad C(t)=\sum_{u=0}^r a_ut^u,\qquad W(t)=\sum_{u=\max\{r+1,\alpha\}}^s a_ut^u,\]
where an empty sum is regarded as zero. Then
\begin{equation}\label{eq:residue-decomposition}
A(t)=(1+t^p+\cdots+t^{qp})C(t)+(1+t^p+\cdots+t^{(q-1)p})W(t).
\end{equation}
Moreover, $h_s=C(1)+W(1)$ and $a=qh_s+C(1)$. If $W(t)\ne0$, then $\alpha\le s$, and the following hold:
\begin{enumerate}[label=\textup{(\roman*)}]
\item if $\alpha<s$, then $C(1)+h_s=h_\alpha+h_{s-1}$;
\item if $\alpha=s$, then $W(t)=wt^s$ for some $w\ne0$ and $h_{s-1}=C(1)$.
\end{enumerate}
\end{lemma}
\begin{proof}
Since $s<p$, we have $d\le R-s$, so $h_i=h_{i-1}$ for $p\le i\le d$. Comparing coefficients in $A(t)\Phi_p(t)$, we obtain $h_i=a_i+a_{i-1}+\cdots+a_{i-p+1}$, and hence
\[0=h_i-h_{i-1}=a_i-a_{i-p}\qquad(p\le i\le d).\]
Applying these equalities successively at the indices $u+p,u+2p,\ldots$, we obtain
\begin{align*}
a_u&=a_{u+p}=\cdots=a_{u+qp}&&(0\le u\le r)\;\;\text{ and}\\
a_u&=a_{u+p}=\cdots=a_{u+(q-1)p}&&(r<u<p).
\end{align*}
Grouping the terms of $A(t)$ by their exponents modulo $p$, we obtain
\[A(t)=(1+t^p+\cdots+t^{qp})C(t)+(1+t^p+\cdots+t^{(q-1)p})\sum_{u=r+1}^{p-1}a_ut^u.\]
We determine which coefficients in the last sum can be nonzero. For $s<u<p$, the constant interval and the coefficient formula imply $0=h_u-h_{u-1}=a_u$.

Next, suppose $r<u<p$ and $u<\alpha$, and set $v=p+r-u$. Then $r<v<p$ and $v>s$. Since $d-u=(q-1)p+v$, palindromicity and the coefficient equalities above imply $a_u=a_{d-u}=a_{(q-1)p+v}=a_v=0$. Hence the last sum equals $W(t)$, proving \eqref{eq:residue-decomposition}. In particular, $h_s=\sum_{u=0}^{p-1}a_u=C(1)+W(1)$. Evaluating the decomposition at $t=1$, we obtain $a=(q+1)C(1)+qW(1)=qh_s+C(1)$.

Suppose $W(t)\ne0$. Since $\alpha=p+r-s>r$, the definition of $W(t)$ implies $\alpha\le s$.
\begin{itemize}
\item[\textup{(i)}] If $\alpha<s$, then $r<\alpha<s<p$. Since $a_u=0$ for $r<u<\alpha$, we have $h_\alpha=C(1)+a_\alpha$. Moreover, palindromicity and the coefficient equalities above imply $a_\alpha=a_{d-\alpha}=a_{(q-1)p+s}=a_s$. Since $s<p$, we also obtain $h_s-h_{s-1}=a_s$. Thus, $h_\alpha-C(1)=a_\alpha=a_s=h_s-h_{s-1}$.
\item[\textup{(ii)}] If $\alpha=s$, then $W(t)=a_st^s$ with $a_s\ne0$. Since $a_u=0$ for $r<u<s$, we have $h_{s-1}=\sum_{u=0}^{s-1}a_u=C(1)$.
\end{itemize}
\end{proof}

\begin{lemma}\label{lem:W-zero}
Under the assumptions of \cref{lem:constant-interval}, let $s$ be the least index with $h_s\le s$ and define $C(t)$ and $W(t)$ as in \cref{lem:residue-decomposition}. Then $W(t)=0$.
\end{lemma}
\begin{proof}
Suppose $W(t)\ne0$. By \cref{lem:constant-interval}, $h_{s-1}=h_s=s$. If $\alpha=s$, then \cref{lem:residue-decomposition}\textup{(ii)} states that $W(t)=wt^s$ with $w\ne0$ and $h_{s-1}=C(1)$. But $w=W(1)=h_s-C(1)=h_s-h_{s-1}=0$, a contradiction. Thus $\alpha<s$, and part \textup{(i)} of the same lemma implies $C(1)=h_\alpha>\alpha$ by minimality of $s$.
Hence, \[a=qs+C(1)>qs+\alpha=p+r+(q-1)s\ge p>a,\]
a contradiction.
\end{proof}

\begin{lemma}\label{lem:replace-prime-general}
Let $a$ be a positive integer, let $p\ge p(a)$ be prime, and let $\ell$ be a positive divisor of $a$. Suppose that $C(t)$ is a product of cyclotomic polynomials with $C(1)=a/\ell$ and $\deg C(t)\le a/\ell-1$. If $C(t)\Psi_{\ell p}(t)$ is CM-realizable and non-type-CI, then so is $C(t)\Psi_{\ell p(a)}(t)$.
\end{lemma}
\begin{proof}
Let $c=C(1)$ and $d=\deg C(t)$. If $d=0$, then $C(t)=1$ and $C(t)\Psi_{\ell p}(t)$ is of type CI, a contradiction. Thus, $d\ge1$. Write $C(t)=\sum_{j=0}^d c_jt^j$ and $u_i=\sum_{j=0}^i c_j$. Since $C(t)$ is monic, we have $u_d=c$ and $u_{d-1}=c-1$. For $m=\ell p$ or $\ell p(a)$, we have $d<m$, and the coefficient sequence of $C(t)\Psi_m(t)$ is
\[u_0,u_1,\ldots,u_{d-1},\underbrace{c,\ldots,c}_{m-d\text{ entries}},u_{d-1},u_{d-2},\ldots,u_0.\]
For $m=\ell p$, this is an $M$-sequence by assumption, and the coefficient of $t^m$ is $u_{d-1}=c-1<m$. By applying \eqref{eq:small-macaulay} repeatedly, we obtain $u_{d-1}\ge\cdots\ge u_0$. Replacing $\ell p$ by $\ell p(a)$ preserves the first $d+1$ coefficients, changes only the length of the constant block, and leaves the subsequent nonincreasing sequence unchanged. Therefore, all Macaulay inequalities remain valid, so $C(t)\Psi_{\ell p(a)}(t)$ is CM-realizable.

For either $q=p$ or $q=p(a)$, we have $d\le a-1<q-1$. If $q\mid n$, then $\ph(n)\ge q-1>d$, so no cyclotomic factor of $C(t)$ has index divisible by $q$. In particular, $b_{\ell q}(C(t))=0$. Multiplication by $\Psi_{\ell q}(t)$ increases this integer from $0$ to $1$ and leaves every other $b_n(C(t))$ unchanged. By \cref{lem:mobius-ci}, $C(t)\Psi_{\ell q}(t)$ is of type CI if and only if $C(t)$ is of type CI. This holds for both choices of $q$, so non-type-CI is preserved as well.
\end{proof}

We are now ready to prove Theorem~\ref{thm:prime-reduction}.
\begin{proof}[Proof of \cref{thm:prime-reduction}]
Let $f(t)$ be as in the theorem. By \cref{lem:p-primary-factor}, we can write $f(t)=A(t)\Phi_p(t)$ with $A(1)=a$. Let $d=\deg A(t)$. We shall find a positive divisor $\ell$ of $a$ and a product $C(t)$ of cyclotomic polynomials such that $f(t)=C(t)\Psi_{\ell p}(t)$ and $\deg C(t)\le C(1)-1=a/\ell-1$.

When $d<p$, we take $\ell=1$ and $C(t)=A(t)$. The degree bound follows from \cref{lem:small-degree-A} with $m=p$.

Now suppose $d\ge p$ and write $d=qp+r$, where $q\ge1$ and $0\le r<p$. By Lemmas~\ref{lem:residue-decomposition} and~\ref{lem:W-zero}, we have
\[A(t)=(1+t^p+\cdots+t^{qp})C(t),\qquad f(t)=C(t)\Psi_{(q+1)p}(t)\;\;\text{ and }\;\; \deg C(t)\le r.\]
The factor $1+t^p+\cdots+t^{qp}=\Psi_{(q+1)p}(t)/\Psi_p(t)$ is a product of cyclotomic polynomials, and so is $A(t)$. Thus, unique factorization in $\QQ[t]$ implies that $C(t)$ is also such a product. Set $\ell=q+1$. Then $A(1)=a=\ell C(1)$, so $\ell\mid a$ and $C(1)=a/\ell$. Since $\deg C(t)\le r<\ell p$, we may apply \cref{lem:small-degree-A} to $C(t)\Psi_{\ell p}(t)$ and conclude that $\deg C(t)\le C(1)-1=a/\ell-1$.

In both cases, \cref{lem:replace-prime-general} applies to $f(t)=C(t)\Psi_{\ell p}(t)$. Therefore, $C(t)\Psi_{\ell p(a)}(t)$ is CM-realizable and non-type-CI, as required.
\end{proof}

\bigskip
\section{Classification of the multiplicities \texorpdfstring{$8k$}{8k}}\label{sec:8k}
In this section, we construct the required families, handle the exceptional prime parameters, and prove \cref{thm:8k-classification}.

\subsection{Existence at every multiplicity \texorpdfstring{$32k$}{32k}}\label{sec:32k}
\begin{proposition}\label{prop:32-seeds}
The following products of cyclotomic polynomials are CM-realizable and non-type-CI:
\begin{align*}
S_0(t)&=\Phi_2(t)^5\Phi_6(t),&S_0(1)&=32,\\
S_1(t)&=\Phi_2(t)^3\Phi_3(t)\Phi_4(t)^2\Phi_{10}(t),&S_1(1)&=96,\\
S_2(t)&=\Phi_2(t)^4\Phi_3(t)\Phi_5(t)\Phi_8(t),&S_2(1)&=480.
\end{align*}
\end{proposition}
\begin{proof}
CM realizability follows by expanding each product and checking Macaulay's inequalities. To prove that the products are non-type-CI, apply \cref{lem:missing-divisor}: $S_0(t)$ contains $\Phi_6(t)$ but not $\Phi_3(t)$; $S_1(t)$ contains $\Phi_{10}(t)$ but not $\Phi_5(t)$; and $S_2(t)$ contains $\Phi_8(t)$ but not $\Phi_4(t)$.
\end{proof}
For illustration, the coefficient sequence of $S_0(t)$ is $(1,4,6,5,5,6,4,1)$. The only increasing steps after index $0$ are from $4$ to $6$ and from $5$ to $6$. These satisfy Macaulay's inequalities because $4^{\langle1\rangle}=10\ge6$ and $5^{\langle4\rangle}=6$, respectively. Hence the sequence is an $M$-sequence. The coefficient sequences of $S_1(t)$ and $S_2(t)$ can be checked in the same way.

\begin{theorem}\label{thm:32k}
For every positive integer $k$, there exists a product $f(t)$ of cyclotomic polynomials which is CM-realizable and non-type-CI and satisfies $f(1)=32k$.
\end{theorem}
\begin{proof}
For a positive integer $n$, write $n=\prod_p p^{\nu_p(n)}$, where all products indexed by $p$ in this proof run over primes. Since $\Phi_p(t)=\Psi_p(t)$ for every prime $p$, multiplication by $\prod_{p\mid n}\Phi_p(t)^{\nu_p(n)}$ preserves CM realizability by the join construction and multiplies the value at $1$ by $n$.

If $3\nmid k$, take $f(t)=S_0(t)\prod_{p\mid k}\Phi_p(t)^{\nu_p(k)}$. Then $f(1)=32k$, and $f(t)$ contains $\Phi_6(t)$ but not $\Phi_3(t)$.

If $3\mid k$ and $5\nmid k$, write $k=3n$ and take $f(t)=S_1(t)\prod_{p\mid n}\Phi_p(t)^{\nu_p(n)}$. Then $f(1)=96n=32k$, and $f(t)$ contains $\Phi_{10}(t)$ but not $\Phi_5(t)$.

Finally, if $15\mid k$, write $k=15n$ and take $f(t)=S_2(t)\prod_{p\mid n}\Phi_p(t)^{\nu_p(n)}$. Then $f(1)=480n=32k$, and $f(t)$ contains $\Phi_8(t)$ but not $\Phi_4(t)$.

In every case, \cref{lem:missing-divisor} guarantees that $f(t)$ is non-type-CI.
\end{proof}

\subsection{The multiplicities \texorpdfstring{$16k$}{16k} for odd \texorpdfstring{$k$}{k}}\label{sec:16k}
\begin{theorem}\label{thm:16k-odd}
Let $k$ be a positive odd integer. Then there exists a product $f(t)$ of cyclotomic polynomials which is CM-realizable and non-type-CI and satisfies $f(1)=16k$ if and only if $k\ne1,3$.
\end{theorem}
\begin{proof}
For $k=1,3$, nonexistence follows from \cref{lem:finite-checks}\textup{(i)}, applied at multiplicities $16,48$. For $k\ge5$, let $F_k(t)=\Phi_2(t)^4\Phi_6(t)\Psi_k(t)$. Then $F_k(1)=16k$. The coefficient sequence of $F_k(t)$ is
\begin{align*}
&(1,4,7,9,12,14,12,9,7,4,1) \;\;\text{ for }\;\; k=5, \\
&(1,4,7,9,12,15,\underbrace{16,\ldots,16}_{k-6\text{ entries}},15,12,9,7,4,1) \;\;\text{ for }\;\;k \geq 7.
\end{align*}
Macaulay's inequalities are verified directly for $k=5$. For $k\ge7$, the coefficients through degree $4$ agree with those for $k=5$, and the only new increases are from $12$ to $15$ and from $15$ to $16$, both of which satisfy Macaulay's inequalities. All subsequent entries are constant or decreasing. Hence, $F_k(t)$ is CM-realizable.

To prove that $F_k(t)$ is non-type-CI, we apply \cref{lem:mobius-ci}. Since $k\ge5$, the observation following that lemma implies $b_3(\Psi_k(t))=0$. Moreover, $b_3(\Phi_6(t))=-1$ and $b_3(\Phi_2(t)^4)=0$. Hence, $b_3(F_k(t))=-1$, which proves that $F_k(t)$ is non-type-CI, as required.
\end{proof}

\subsection{The multiplicities \texorpdfstring{$8pq$}{8pq} for odd \texorpdfstring{$p,q$}{p,q}}
\begin{proposition}\label{prop:Fpq}
Let $7\le p\le q$ be odd primes. The product of cyclotomic polynomials defined by $F_{p,q}(t)=\Phi_2(t)^3\Phi_{10}(t)\Phi_p(t)\Phi_q(t)$ is CM-realizable and non-type-CI, and satisfies $F_{p,q}(1)=8pq$.
\end{proposition}
\begin{proof}
The value at $1$ is immediate. Since $p,q\ge7$, the product contains $\Phi_{10}(t)$ but not $\Phi_5(t)$, so it is non-type-CI by \cref{lem:missing-divisor}.

For CM realizability, use $\Phi_2(t)\Phi_{10}(t)=1+t^5$ to write $F_{p,q}(t)=(1+t)^2(1+t^5)\Psi_p(t)\Psi_q(t)$. If we write $\Psi_p(t)\Psi_q(t)=\sum_{j\in\ZZ}c_jt^j$, then
\[c_j=\begin{cases}j+1,&0\le j\le p-1,\\p,&p\le j\le q-1,\\p+q-1-j,&q\le j\le p+q-2,\\0,&\text{otherwise}.\end{cases}\]
Writing $F_{p,q}(t)=\sum_i h_it^i$, we obtain \[h_i=c_i+2c_{i-1}+c_{i-2}+c_{i-5}+2c_{i-6}+c_{i-7}.\]
All coefficients are nonnegative, and $(h_0,\ldots,h_6)=(1,4,8,12,16,21,28)$. The inequalities for $1\le i\le5$ follow from
\[4^{\langle1\rangle}=10,\quad8^{\langle2\rangle}=13,\quad12^{\langle3\rangle}=17,\quad16^{\langle4\rangle}=22,\quad21^{\langle5\rangle}=28.\]
The inequalities at the next two positions follow from $28^{\langle6\rangle}=36$ and the following values:
\[\begin{array}{c|ccc}
(p,q)&(7,7)&p=7<q&p\ge11\\\hline
h_7&34&35&36\\
h_7^{\langle7\rangle}&40&42&45\\
h_8&36&40&44
\end{array}\]
For $8\le i\le p-2$, direct substitution shows that $h_i=8i-20=8i-\binom{7}{2}+1$ and $h_{i+1}-h_i=8$, so we can apply \cref{lem:macaulay-increment}.

Consider the remaining indices $p-1\le i\le p+5$ with $i\ge8$, and let $d_i=h_{i+1}-h_i$. When $d_i>0$, the same coefficient formula implies the bounds
\[\begin{array}{c|cc}
i&\text{lower bound for }h_i&\text{upper bound for }d_i\\\hline
p-1&8p-28&7\\
p&8p-22&5\\
p+1&8p-16&4\\
p+2&8p-13&4\\
p+3&8p-8&4\\
p+4&8p-4&3\\
p+5&8p-1&1
\end{array}\]
For $1\le d\le i$, the expression $di-\binom{d-1}{2}+1$ is increasing in $d$. In every displayed row with $i\ge8$, its value at the upper bound for $d_i$ does not exceed the lower bound for $h_i$, since $p\ge7$. Thus \cref{lem:macaulay-increment} applies to each positive step; the indices below $8$ have already been checked.

Finally, $c_{j+1}-c_j\le0$ for $j\ge p-1$. If $i\ge p+6$, all indices $i,i-1,i-2,i-5,i-6,i-7$ are at least $p-1$, so $h_{i+1}\le h_i$. Hence all remaining Macaulay inequalities hold, and the coefficient sequence is an $M$-sequence.
\end{proof}

\begin{lemma}\label{lem:P5-P9}
Let $P_5(t)=\Phi_2(t)^2\Phi_5(t)\Phi_8(t)$ and $P_9(t)=\Phi_2(t)^2\Phi_3(t)^2\Phi_8(t)$. These products of cyclotomic polynomials are CM-realizable and non-type-CI, with $P_5(1)=40$ and $P_9(1)=72$.
\end{lemma}
\begin{proof}
The coefficient sequences of both products satisfy Macaulay's inequalities by direct computation. Each product contains $\Phi_8(t)$ but not $\Phi_4(t)$, so both are non-type-CI by \cref{lem:missing-divisor}.
\end{proof}

\begin{theorem}\label{thm:8pq}
Let $p,q\ge5$ be odd integers. Then there exists a product $f(t)$ of cyclotomic polynomials which is CM-realizable and non-type-CI and satisfies $f(1)=8pq$.
\end{theorem}
\begin{proof}
If $5\mid pq$, take $f(t)=P_5(t)\Psi_{pq/5}(t)$. If $5\nmid pq$ and $9\mid pq$, take $f(t)=P_9(t)\Psi_{pq/9}(t)$. Both products are CM-realizable and have value $8pq$ at $1$. In each case, the added $\Psi$-factor has odd index, so it does not introduce $\Phi_4(t)$, while $\Phi_8(t)$ remains. Hence $f(t)$ is non-type-CI in both cases.

Suppose $\gcd(pq,15)=1$. Choose prime factors $7\le r\le s$ of $pq$, counted with multiplicity, so that $rs\mid pq$ and take $f(t)=F_{r,s}(t)\Psi_{pq/rs}(t)$. The product is CM-realizable, has value $8pq$, and contains $\Phi_{10}(t)$ but not $\Phi_5(t)$.

It remains to consider the case $5\nmid pq$, $3\mid pq$, and $9\nmid pq$. Exactly one of $p,q$ is divisible by $3$. After removing this factor of $3$, the integer $pq/3$ still has at least two prime factors, counted with multiplicity, each at least $7$. Choose primes $7\le r\le s$ with $rs\mid pq/3$ and take $f(t)=F_{r,s}(t)\Psi_{pq/3rs}(t)\Psi_3(t)$. Then $f(t)$ is CM-realizable, satisfies $f(1)=8pq$, and contains $\Phi_{10}(t)$ but not $\Phi_5(t)$.
\end{proof}
\begin{remark}\label{rem:4pq-nonexistence}
By \cref{thm:8pq}, for every pair of odd primes $p,q\ge5$, there exists a product of cyclotomic polynomials which is CM-realizable and non-type-CI and has multiplicity $8pq$. The analogous assertion with $4$ in place of $8$ is false. In fact, \cref{thm:prime-reduction} and \cref{lem:finite-checks}\textup{(ii)}, applied with $a=20,28$, show that no such product has multiplicity $20p$ for a prime $p\ge23$, or $28p$ for a prime $p\ge29$.
\end{remark}

\subsection{The exceptional multiplicities \texorpdfstring{$24p$}{24p}}\label{subsec:24p}

We now consider $p=17,19,23$, which are not covered by \cref{thm:prime-reduction} with $a=24$.

\begin{lemma}\label{lem:24p-phi-p-reduction}
Let $p\in\{17,19,23\}$, and let $A(t)$ be a product of cyclotomic polynomials with $A(1)=24$. If $A(t)\Phi_p(t)$ is CM-realizable, then there exist a positive divisor $\ell$ of $24$ and a product $C(t)$ of cyclotomic polynomials such that
\[A(t)\Phi_p(t)=C(t)\Psi_{\ell p}(t),\qquad C(1)=\frac{24}{\ell},\qquad\deg C(t)\le\frac{24}{\ell}-1.\]
\end{lemma}
\begin{proof}
Write $A(t)=\sum_{i=0}^d a_it^i$ and $A(t)\Phi_p(t)=\sum_{i=0}^R h_it^i$, where $R=d+p-1$. If $d\le23$, take $\ell=1$ and $C(t)=A(t)$. Henceforth, assume $d\ge24$.

Let $s$ be the smallest index with $h_s\le s$. If $2s>R$, then minimality of $s$ implies $h_i\ge i+1$ for $0\le i\le\lfloor R/2\rfloor$. Since $A(t)\Phi_p(t)$ is palindromic, $R\ge p+23$, and $p+23$ is even, we have
\[24p=\sum_{i=0}^R h_i\ge\left(\left\lfloor\frac R2\right\rfloor+1\right)^2\ge\left(\frac{p+25}{2}\right)^2>24p,\]
a contradiction for each $p\in\{17,19,23\}$. Thus, $2s\le R$. Macaulay's inequalities and palindromicity now imply $h_s=\cdots=h_{R-s}$. The inequality at index $R-s$, together with palindromicity and minimality of $s$, also implies $h_{s-1}=h_s=s$. Also, $R=d+p-1 \leq d+22 <2d$, so $s<d$.

We claim that $s<p$. Suppose otherwise, and put $P_j=\sum_{i=0}^j a_i$, with $P_j=0$ for $j<0$. For $0\le j\le d$, the coefficient formula for $h_j$ implies $P_j=h_j+P_{j-p}$. Since $h_j\ge0$, induction on $j$ shows that $P_j\ge h_j\ge0$. Splitting the coefficients of $A(t)$ into three consecutive blocks and using palindromicity, we obtain
\[24=A(1)=\sum_{i=0}^{s-p-1}a_i+\sum_{i=s-p}^{s-1}a_i+\sum_{i=s}^d a_i=P_{s-p-1}+h_{s-1}+P_{d-s}\ge s+h_{d-s}.\]
Since $0\le d-s\le R-s$, minimality of $s$ and the constant interval imply $h_{d-s}\ge\min\{d-s+1,s\}$. Therefore, $24\ge s+h_{d-s}\ge\min\{d+1,2s\}>24$, since $d\ge24$ and $s\ge p\ge17$. This contradiction proves $s<p$.

Write $d=qp+r$ with $q\ge1$ and $0\le r<p$, and apply \cref{lem:residue-decomposition}. With the notation of that lemma, we have
\[A(t)=(1+t^p+\cdots+t^{qp})C(t)+(1+t^p+\cdots+t^{(q-1)p})W(t) \;\;\text{ and }\;\; 24=qs+C(1).\]
Suppose $W(t)\ne0$. If $\alpha=s$, then part \textup{(ii)} of that lemma implies $W(t)=wt^s$ for some $w\ne0$ and $h_{s-1}=C(1)$. Since $h_s=h_{s-1}$ and $h_s=C(1)+W(1)$, we obtain $W(1)=w=0$, a contradiction. Hence $\alpha<s$, and part \textup{(i)}, together with $h_s=h_{s-1}$ and minimality of $s$, implies $C(1)=h_\alpha>\alpha$. Consequently,
\[24=qs+C(1)>qs+\alpha=p+r+(q-1)s.\]
For $q=1$, the right-hand side equals $d\ge24$, a contradiction. For $q\ge2$, the inequality $\alpha<s$ implies $2s>p+r$, so the right-hand side exceeds $3(p+r)/2\ge3\cdot17/2>24$, again a contradiction. Therefore, $W(t)=0$.

Hence, $A(t)\Phi_p(t)=C(t)\Psi_{(q+1)p}(t)$. The factor $1+t^p+\cdots+t^{qp}$ is a product of cyclotomic polynomials, so unique factorization shows that $C(t)$ is also such a product. Set $\ell=q+1$. Then $24=A(1)=\ell C(1)$, and $\deg C(t)\le r<\ell p$. Finally, \cref{lem:small-degree-A} can be applied to $C(t)\Psi_{\ell p}(t)$, proving $\deg C(t)\le C(1)-1=24/\ell-1$.
\end{proof}

\begin{lemma}\label{lem:24p-phi-p2}
Let $p\in\{17,19,23\}$ and let $A(t)$ be a product of cyclotomic polynomials with $A(1)=24$. Then $A(t)\Phi_{p^2}(t)$ is not CM-realizable.
\end{lemma}
\begin{proof}
Suppose that $f(t)=A(t)\Phi_{p^2}(t)=\sum_{i=0}^R h_it^i$ is CM-realizable. We have $R\ge p(p-1)$, $h_0=h_R=1$ and $h_{R-1}=h_1$. If $h_1\ge2$, then $h_i\ge2$ for every $1\le i\le R-1$: once an entry is at most $1$, Macaulay's inequalities force all subsequent entries to be at most $1$, contradicting $h_{R-1}=h_1\ge2$. Hence,
\[24p=f(1)\ge2R\ge2p(p-1) \geq 2p \cdot 16 >24p,\]
a contradiction. Thus $h_1\le1$. Since $h_R=1$, Macaulay's inequalities force $h_i=1$ for all $0\le i\le R$. Hence $f(t)=\Psi_{R+1}(t)=\Psi_{24p}(t)$. However, $\Phi_{p^2}(t)$ divides $f(t)$, whereas $p^2\nmid24p$, a contradiction.
\end{proof}

\begin{proposition}\label{prop:24p-small-primes}
For $p\in\{17,19,23\}$, there is no product $f(t)$ of cyclotomic polynomials which is CM-realizable and non-type-CI and satisfies $f(1)=24p$.
\end{proposition}
\begin{proof}
Suppose such a product exists. Since $p\nmid24$, exactly one cyclotomic factor has index a power of $p$, counted with multiplicity; write it as $\Phi_{p^e}(t)$. If $e\ge3$, then $\deg f(t)\ge p^2(p-1)>24p-1$, a contradiction to \cref{lem:degree-bound}. The case $e=2$ is excluded by \cref{lem:24p-phi-p2}. Thus, $e=1$, and $f(t)=A(t)\Phi_p(t)$ with $A(1)=24$. Lemma~\ref{lem:24p-phi-p-reduction} reduces this case to the finite set of candidates excluded by \cref{lem:finite-checks}\textup{(iii)}.
\end{proof}

\subsection{Three additional explicit polynomials}
\begin{proposition}\label{prop:Q-seeds}
The following products of cyclotomic polynomials are CM-realizable and non-type-CI:
\begin{align*}
Q_7(t)&=\Psi_2(t)\Psi_3(t)\Psi_4(t)\Psi_7(t)\Phi_{12}(t),&Q_7(1)&=168,\\
Q_{11}(t)&=\Psi_3(t)\Psi_4(t)\Psi_{11}(t)\Phi_4(t)\Phi_{20}(t),&Q_{11}(1)&=264,\\
Q_{13}(t)&=\Psi_3(t)\Psi_8(t)\Psi_{13}(t)\Phi_{24}(t),&Q_{13}(1)&=312.
\end{align*}
\end{proposition}
\begin{proof}
CM realizability follows by expanding the products and checking Macaulay's inequalities. To prove that the products are non-type-CI, apply \cref{lem:missing-divisor}: $Q_7(t)$ contains $\Phi_{12}(t)$ but not $\Phi_6(t)$; $Q_{11}(t)$ contains $\Phi_{20}(t)$ but not $\Phi_{10}(t)$; and $Q_{13}(t)$ contains $\Phi_{24}(t)$ but not $\Phi_{12}(t)$.
\end{proof}

\subsection{Completing the classification}
We are now ready to prove \cref{thm:8k-classification}.
\begin{proof}[Proof of \cref{thm:8k-classification}]
Assume $k\notin\Nc$. If $k$ is even and $4\mid k$, then we apply \cref{thm:32k} to $8k=32(k/4)$. If $k\equiv2\pmod4$, write $k=2\ell$ with $\ell$ odd. Since $k\notin\{2,6\}$, we have $\ell\ne1,3$, so \cref{thm:16k-odd} provides the required example at multiplicity $8k=16\ell$.

Let $k$ be odd. If $5\mid k$, use $P_5(t)\Psi_{k/5}(t)$; if $9\mid k$, use $P_9(t)\Psi_{k/9}(t)$, where $P_5(t)$ and $P_9(t)$ are defined in \cref{lem:P5-P9}. In each case, the product is CM-realizable, has multiplicity $8k$, and contains $\Phi_8(t)$ but not $\Phi_4(t)$. It remains to consider odd $k$ with $5\nmid k$ and $9\nmid k$.
\begin{itemize}
\item If $3\nmid k$, then $k\notin\Nc$ implies that $k$ is composite. Thus, $k=uv$ for odd $u,v\ge5$, and \cref{thm:8pq} guarantees the desired polynomial.
\item If $3\mid k$, write $k=3m$ with $\gcd(m,15)=1$. When $m$ is composite, write $m=uv$ with odd $u,v\ge7$ and apply \cref{thm:8pq} to $k=(3u)v$. When $m$ is prime, the assumption $k\notin\Nc$ forces $m\in\{7,11,13\}$, and the polynomials in \cref{prop:Q-seeds} provide the required examples. The remaining case $m=1$ is excluded by $3\in\Nc$.
\end{itemize}

Conversely, let $k\in\Nc$. For $k=1,2,3,6$, the multiplicities $8k=8,16,24,48$ are excluded by \cref{lem:finite-checks}\textup{(i)}. If $k=p\ge7$ is prime, the same finite verification excludes multiplicity $8p=56$ when $p=7$. For $p\ge11=p(8)$, \cref{thm:prime-reduction} with $a=8$ and \cref{lem:finite-checks}\textup{(ii)} exclude multiplicity $8p$.

Finally, suppose $k=3p$ for a prime $p\ge17$. The cases $p=17,19,23$ are excluded by \cref{prop:24p-small-primes}. For $p\ge29=p(24)$, apply \cref{thm:prime-reduction} with $a=24$ and \cref{lem:finite-checks}\textup{(ii)}.
\end{proof}

\bigskip
\section{Further consequences}\label{sec:further}
In this section, we give partial existence and nonexistence results for $4k$ and $2k$ with $k$ odd, determine the multiplicities below $100$, and discuss the smallest odd multiplicity and Gorenstein realizability.

\subsection{Existence families for \texorpdfstring{$4k$}{4k} with \texorpdfstring{$k$}{k} odd}
The following products of cyclotomic polynomials are CM-realizable and non-type-CI:
\begin{center}
\begin{tabular}{@{}clc@{}}
\toprule
$b$ & $F_b(t)$ & $F_b(1)$\\
\midrule
$15$ & $\Phi_2(t)\Phi_3(t)\Phi_5(t)\Phi_8(t)$ & $4\cdot15$\\
$21$ & $\Phi_2(t)\Phi_3(t)\Phi_4(t)\Phi_7(t)\Phi_{12}(t)$ & $4\cdot21$\\
$27$ & $\Phi_2(t)\Phi_3(t)^3\Phi_8(t)$ & $4\cdot27$\\
$35$ & $\Phi_2(t)^2\Phi_5(t)\Phi_7(t)\Phi_{12}(t)$ & $4\cdot35$\\
\bottomrule
\end{tabular}
\end{center}
Their coefficient sequences satisfy Macaulay's inequalities by direct computation. Each product contains either $\Phi_8(t)$ without $\Phi_4(t)$, or $\Phi_{12}(t)$ without $\Phi_6(t)$.

If a positive odd integer $k$ is divisible by one of $15,21,27,35$, then the corresponding product $F_b(t)\Psi_{k/b}(t)$ is CM-realizable, has value $4k$ at $1$, and retains the same missing cyclotomic factor. On the other hand, no product of cyclotomic polynomials which is CM-realizable and non-type-CI has multiplicity $4p$ for an odd prime $p$. The case $p=3$ is excluded by \cref{lem:finite-checks}\textup{(i)} at multiplicity $12$; for $p\ge5$, apply \cref{thm:prime-reduction} with $a=4$ and \cref{lem:finite-checks}\textup{(ii)}.

These results are still far from a complete classification of the multiplicities $4k$ with $k$ odd. In fact, even when $k$ is the product of two odd primes, both existence and nonexistence occur. For example, our constructions cover $k=15,35$, whereas \cref{rem:4pq-nonexistence} shows nonexistence for $k=5p$ for every prime $p\ge23$. Moreover, although \cref{thm:prime-reduction} reduces the problem to finitely many cases when one prime factor is fixed and the other is sufficiently large, these finite cases depend on the fixed prime. Thus, a complete classification appears to require a more uniform approach.

\subsection{Existence families for \texorpdfstring{$2k$}{2k} with \texorpdfstring{$k$}{k} odd}
Similarly, the following products of cyclotomic polynomials are CM-realizable and non-type-CI:
\begin{center}
\begin{tabular}{@{}clc@{}}
\toprule
$b$ & $G_b(t)$ & $G_b(1)$\\
\midrule
$27$ & $\Phi_3(t)^3\Phi_8(t)$ & $2\cdot27$\\
$35$ & $\Phi_2(t)\Phi_5(t)\Phi_7(t)\Phi_{12}(t)$ & $2\cdot35$\\
$45$ & $\Phi_3(t)^2\Phi_5(t)\Phi_8(t)$ & $2\cdot45$\\
\bottomrule
\end{tabular}
\end{center}
Macaulay's inequalities can again be checked by expanding the products, and the same missing-divisor argument applies. Therefore, whenever a positive odd integer $k$ is divisible by $27$, $35$, or $45$, multiplication of the corresponding $G_b(t)$ by $\Psi_{k/b}(t)$ produces an example of multiplicity $2k$.

For an odd prime $p$, no product of cyclotomic polynomials which is CM-realizable and non-type-CI has multiplicity $2p$ or $6p$. For $2p$, apply \cref{thm:prime-reduction} with $a=2$ and \cref{lem:finite-checks}\textup{(ii)}. For $6p$, the same argument with $a=6$ handles $p\ge7$, and \cref{lem:finite-checks}\textup{(i)} excludes $p=3,5$ at multiplicities $18,30$. These statements do not constitute complete classifications at multiplicities $2k$ or $4k$; the families and exclusions above do not settle all odd values of $k$.

\subsection{Small multiplicities and the smallest odd value}
By \cref{lem:finite-checks}\textup{(i)}, the multiplicities $1\le n\le99$ admitting a product of cyclotomic polynomials which is CM-realizable and non-type-CI are precisely
\[\{32,40,54,60,64,70,72,80,84,90,96\}.\]
In particular, all such multiplicities below $100$ are even.

On the other hand, the smallest odd multiplicity of a product of cyclotomic polynomials which is CM-realizable and non-type-CI is $105$.
To construct an example at this multiplicity, set \[f(t)=\Phi_3(t)\Phi_5(t)\Phi_7(t)\Phi_{12}(t).\] Its value at $1$ is $105$,
and its coefficient sequence satisfies Macaulay's inequalities. The product contains $\Phi_{12}(t)$ but not $\Phi_6(t)$, so it is non-type-CI.

The finite enumeration excludes all smaller odd multiplicities through $99$. Prime multiplicities are excluded by \cref{thm:prime-reduction} with $a=1$: the only possible choice is $\ell=1$ and $C(t)=1$, and $C(t)\Psi_2(t)$ is of type CI. In particular, $101$ and $103$ are excluded, completing the proof of minimality.

\subsection{Gorenstein realizability}\label{subsec:gorenstein}
In this subsection, the ground field is $\QQ$. A polynomial $f(t)$ is called \emph{Gorenstein-realizable} (resp. \emph{Gorenstein$^*$-realizable}) if it is the $h$-polynomial of a Gorenstein (resp. Gorenstein$^*$) simplicial complex over $\QQ$. Thus, these terms refer specifically to Stanley--Reisner realizations, not to realizations by arbitrary standard graded Gorenstein rings.

Recall that a simplicial complex is Gorenstein$^*$ over $\QQ$ if its Stanley--Reisner ring is Gorenstein with $a$-invariant $0$. Every Gorenstein complex is the join of its Gorenstein$^*$ core with a simplex, and this join does not change the $h$-polynomial; see \cite[Section~5.6]{BrunsHerzog}. Hence, Gorenstein realizability and Gorenstein$^*$ realizability are equivalent at the level of $h$-polynomials.

For a nonconstant product $f(t)=\sum_{i=0}^d h_it^i$ of cyclotomic polynomials, define $g_0=1$ and $g_i=h_i-h_{i-1}$ for $1\le i\le\lfloor d/2\rfloor$. The sequence $(h_0,\ldots,h_d)$ is palindromic. Adiprasito's Lefschetz theorem for rational homology spheres \cite[Theorem~9.1]{Adiprasito}, together with the sufficiency part of the $g$-theorem \cite{BilleraLee,Stanley1980}, implies that
\begin{equation}\label{eq:gorenstein-criterion}
f(t)\text{ is Gorenstein-realizable}\quad\Longleftrightarrow\quad(g_0,\ldots,g_{\lfloor d/2\rfloor})\text{ is an }M\text{-sequence}.
\end{equation}

Consider the following products of cyclotomic polynomials:
\begin{align*}
T_0(t)&=\Phi_2(t)^4\Phi_4(t)\Phi_6(t),&T_0(1)&=32,\\*
T_1(t)&=\Phi_2(t)^4\Phi_3(t)\Phi_4(t)\Phi_6(t)\Phi_{10}(t),&T_1(1)&=96,\\*
T_2(t)&=\Phi_2(t)^4\Phi_3(t)\Phi_5(t)\Phi_8(t),&T_2(1)&=480.
\end{align*}
Their $g$-vectors are $(1,2,1,1,1)$, $(1,2,2,2,1,1,1,0)$, and $(1,5,11,14,12,9,7,3)$, respectively. Each of them is an $M$-sequence, so all three products are Gorenstein-realizable. They are non-type-CI because $T_0(t)$ contains $\Phi_6(t)$ without $\Phi_3(t)$, $T_1(t)$ contains $\Phi_{10}(t)$ without $\Phi_5(t)$, and $T_2(t)$ contains $\Phi_8(t)$ without $\Phi_4(t)$.

\begin{proposition}\label{prop:gorenstein-32k}
For every positive integer $k$, there exists a product $f(t)$ of cyclotomic polynomials which is Gorenstein-realizable and non-type-CI and satisfies $f(1)=32k$.
\end{proposition}
\begin{proof}
Joins with the boundary complexes of simplices preserve the Gorenstein$^*$ property. Thus, multiplication by $\Phi_p(t)=\Psi_p(t)$ preserves Gorenstein realizability for every prime $p$. Repeat the proof of \cref{thm:32k} by replacing each $S_i(t)$ by $T_i(t)$ for $i=0,1,2$. The same missing cyclotomic factors show that the resulting products are non-type-CI.
\end{proof}

Two further examples are $\Phi_2(t)\Phi_3(t)\Phi_5(t)\Phi_8(t)$, of multiplicity $60$, and $\Phi_3(t)^2\Phi_5(t)\Phi_8(t)$, of multiplicity $90$. Their $g$-vectors are $(1,2,2,1,1,1)$ and $(1,2,3,2,2,1,1)$, respectively. Both are $M$-sequences, and both products contain $\Phi_8(t)$ without $\Phi_4(t)$.

By \cref{lem:finite-checks}\textup{(iv)}, the multiplicities $1\le n\le99$ admitting a product of cyclotomic polynomials which is Gorenstein-realizable and non-type-CI are precisely $32,60,64,90,96$. Thus, the multiplicities $40,54,70,72,80,84$ admit products of cyclotomic polynomials which are CM-realizable and non-type-CI, but admit no products which are Gorenstein-realizable and non-type-CI. In particular, $\Phi_2(t)^4\Phi_5(t)\Phi_6(t)$ from \cref{thm:16k-odd} is a product of cyclotomic polynomials which is CM-realizable and non-type-CI and has multiplicity $80$, whereas no product of cyclotomic polynomials of that multiplicity is both Gorenstein-realizable and non-type-CI in the sense of this subsection.

\bigskip
\appendix
\section{Finite computations}\label{sec:computation}
We collect the finite computations used in the proofs and describe exhaustive procedures for reproducing them. The coefficient calculations for individual explicit examples require only expansion and Macaulay's inequalities.

\begin{lemma}\label{lem:finite-checks}
The following assertions hold.
\begin{enumerate}[label=\textup{(\roman*)}]
\item For $1\le n\le99$, there exists a product $f(t)$ of cyclotomic polynomials which is CM-realizable and non-type-CI and satisfies $f(1)=n$ if and only if
\[n\in\{32,40,54,60,64,70,72,80,84,90,96\}.\]
\item For $a\in\{2,4,6,8,20,24,28\}$, there is no pair consisting of a positive divisor $\ell$ of $a$ and a product $C(t)$ of cyclotomic polynomials satisfying $C(1)=a/\ell$ and $\deg C(t)\le a/\ell-1$ such that $C(t)\Psi_{\ell p(a)}(t)$ is CM-realizable and non-type-CI.
\item For $p\in\{17,19,23\}$, there is no pair consisting of a positive divisor $\ell$ of $24$ and a product $C(t)$ of cyclotomic polynomials satisfying $C(1)=24/\ell$ and $\deg C(t)\le24/\ell-1$ such that $C(t)\Psi_{\ell p}(t)$ is CM-realizable and non-type-CI.
\item For $1\le n\le99$, there exists a product $f(t)$ of cyclotomic polynomials which is Gorenstein-realizable and non-type-CI and satisfies $f(1)=n$ if and only if $n\in\{32,60,64,90,96\}$.
\end{enumerate}
\end{lemma}
\begin{proof}[Finite verification]
We use two exhaustive searches: one over coefficient sequences for \textup{(i)} and \textup{(iv)}, and one over cyclotomic factorizations for \textup{(ii)} and \textup{(iii)}.

\noindent\textup{(i)} Fix $1\le n\le99$. Every product of cyclotomic polynomials which is CM-realizable and has value $n$ at $1$ has a palindromic coefficient sequence of positive integers $(h_0,\ldots,h_s)$ with $h_0=h_s=1$ and $\sum h_i=n$, so $s\le n-1$. Enumerate these finite sequences, retain those satisfying Macaulay's inequalities, and factor the corresponding polynomials over $\QQ$. Apply \cref{lem:mobius-ci} whenever all irreducible factors are cyclotomic. The values for which at least one product is non-type-CI are exactly those in \textup{(i)}.

\noindent\textup{(ii)} Fix $a\in\{2,4,6,8,20,24,28\}$ and a positive divisor $\ell$ of $a$, and put $c=a/\ell$. Enumerate all products $C(t)=\prod_{m>1}\Phi_m(t)^{e_m}$ satisfying
\[\sum_{m>1}e_m\ph(m)\le c-1,\qquad \prod_{m>1}\Phi_m(1)^{e_m}=c,\qquad e_m\in\ZZ_{\ge0}.\]
For $c=1$, only $C(t)=1$ occurs. Otherwise, $\ph(m)\le c-1$ for every occurring factor, and $m\le2\ph(m)^2\le2(c-1)^2$; the degree bound also limits the exponents. Thus, the enumeration is finite and exhaustive. For each candidate, check $C(t)\Psi_{\ell p(a)}(t)$ against Macaulay's inequalities and \cref{lem:mobius-ci}. No candidate meets both required conditions.

\noindent\textup{(iii)} Repeat the search in \textup{(ii)} with $a=24$, using $C(t)\Psi_{\ell p}(t)$ for $p=17,19,23$. Again no candidate meets both conditions.

\noindent\textup{(iv)} Apply the $g$-vector criterion \eqref{eq:gorenstein-criterion} to every product of cyclotomic polynomials retained in \textup{(i)} which is CM-realizable and non-type-CI. This yields precisely the values in \textup{(iv)}.
\end{proof}

\end{document}